\documentclass[12pt]{article}

\usepackage{amstext    }
\usepackage{amsthm    }
\usepackage{a4}
\usepackage[mathscr]{eucal}

\usepackage{amsmath}
\usepackage{amssymb}
\usepackage{amscd}

\numberwithin{equation}{section}

\newtheorem{theorem}{Theorem}[section]

\newtheorem{proposition}[theorem]{Proposition}
\newtheorem{corollary}[theorem]{Corollary}
\newtheorem{lemma}[theorem]{Lemma}

\newtheorem{conjecture}[theorem]{Conjecture}

\newcommand{\cali}[1]{\mathscr{#1}}

\newcommand{\dist}{{\rm \ \! dist}}

\newcommand{\ddc}{dd^c}

\newcommand{\Cc}{\cali{C}}

\newcommand{\Ec}{\cali{E}}
\newcommand{\Fc}{\cali{F}}

\newcommand{\Hc}{\cali{H}}

\newcommand{\FS}{{\rm FS}}

\newcommand{\C}{\mathbb{C}}

\newcommand{\R}{\mathbb{R}}
\renewcommand{\P}{\mathbb{P}}

\title{Equidistribution of varieties for endomorphisms of projective spaces}
\author{Tien-Cuong Dinh and Nessim Sibony}
\date{{\it Dedicated to Professor Ha Huy Khoai}}

\begin{document}
\maketitle

\begin{abstract}
Let $f$ be a non-invertible holomorphic endomorphism of the complex projective space $\P^k$ and $f^n$ its
iterate of order $n$. 
Let $V$ be an algebraic subvariety of $\P^k$ which is generic in the Zariski sense. We give here a survey on the asymptotic equidistribution of the sequence $f^{-n}(V)$ when $n$ goes to infinity.
\end{abstract}




\section{Introduction}

Let $\P^k$ denote the complex projective space of dimension $k$. Consider an endomorphism
$f:\P^k\rightarrow\P^k$ which is holomorphic and non-constant. Such a map is 
always induced by a polynomial map $F=(F_0,\ldots,F_k)$ from $\C^{k+1}$ to $\C^{k+1}$ where the $F_i$ are homogeneous polynomials of the same degree such that $F^{-1}(0)=\{0\}$. Indeed, if $\pi:\C^{k+1}\setminus\{0\}\to\P^k$ is the canonical projection, the map $f$ is defined by the relation 
$f\circ\pi=\pi\circ F$. 
We refer to \cite{DinhSibony4,Fornaess,Sibony} for the basic properties of these endomorphisms.

From now on, we assume that the algebraic degree of $f$, i.e. the common degree $d$ of the $F_i$, is at least 2. Otherwise, $f$ corresponds to an invertible matrix and its dynamics is easy to study. 
The parameter space for these endomorphisms with a given algebraic degree $d$ is a Zariski open set  of a projective space $\P^N$ that we denote by $\Hc_d(\P^k)$. 
Using the B\'ezout theorem, it is not difficult to see that
an endomorphism $f$ as above defines a ramified covering of
degree $d^k$ over $\P^k$.  In other words, $f^{-1}(a)$ contains exactly $d^k$ points counted with multiplicity. 

Let $\omega_\FS$ denote the Fubini-Study
form on $\P^k$ normalized so that $\omega_\FS^k$ is a probability
measure. Let $f^n:=f\circ \cdots\circ f$ ($n$ times) be the iterate of
order $n$ of $f$. It is well-known that the sequence of probability measures 
$d^{-kn} (f^n)^*(\omega_\FS^k)$
converges to a probability measure $\mu$ which is totally invariant:
$d^{-k}f^*(\mu)=f_*(\mu)=\mu$. It is called {\it the equilibrium measure} or {\it the Green measure} of $f$, see e.g. \cite{DinhSibony4}.

Consider a point  $a$ in $\P^k$. 
We are interested in the asymptotic distribution of the fibers
$f^{-n}(a)$ of $a$ when $n$ goes to infinity. 
We will survey results on the equidistribution of these sets.
The proof for the main results in this section was given in \cite{DinhSibony5}. We will sketch it in Section \ref{section_proof} with some simplification of the arguments. Equidistribution for higher dimension subvarieties will be discussed in Section \ref{section_varieties}. We also refer to Yuan \cite{Yuan} for analogous equidistribution problems in number theory.

We have the following result which was proved in \cite{DinhSibony1}, see also \cite{Dinh}.

\begin{theorem} \label{th_exceptional}
Let $f$ be a non-invertible holomorphic endomorphism of $\P^k$. Then, there is a finite number of algebraic 
subsets  $E\subset \P^k$ which are totally invariant, i.e. $f^{-1}(E)=E=f(E)$. In particular, there is a maximal proper algebraic subset $\Ec$, possibly empty, which is totally invariant.
\end{theorem}
 
Note that $E$ is totally invariant if and only if $f^{-1}(E)\subset E$. We do not assume here that $E$ is irreducible nor of pure dimension. The set $\Ec$ is in fact the union of all totally invariant proper algebraic sets of $f$. These sets are a posteriori of bounded degree and we can construct them explicitly. However, they are far from being understood. 
The following folklore conjecture is still open in dimension $k\geq 3$, see \cite{FornaessSibony3,CerveauLinsNeto} for the dimension 2 case.

\begin{conjecture}
Any totally invariant algebraic subset for a map $f$ as above is a union of linear projective subspaces.
\end{conjecture} 

One also expects that the degrees of these totally invariant sets are bounded by a constant which depends only on $k$. This is known for the case where the codimension of $E$ is 1 or 2 or in some others situations, see \cite{AmerikCampana,DinhSibony1, DinhSibony4, FornaessSibony3}. In dimension 1, $\Ec$ contains 0,1 or 2 points, e.g. if $f(z)=z^{\pm d}$ then $\Ec=\{0,\infty\}$.

Denote by $\delta_a$ the Dirac mass at $a$ and
$\mu_n^a:=d^{-kn}(f^n)^*(\delta_a)$ the probability measure
which is equidistributed on the fiber $f^{-n}(a)$. The points in
$f^{-n}(a)$ are counted with multiplicity.
Here is the first main result \cite{DinhSibony5}.

\begin{theorem} \label{th_equi_2} 
Let $f$, $\mu$ and $\Ec$ be as above. 
There is a constant $\lambda>1$ such that
if $a$ is a point out of $\Ec$, then $\mu_n^a$ converges to $\mu$ exponentially fast, that is, 
if  $\varphi$ is a $\Cc^\alpha$ function on $\P^k$
with $0<\alpha \leq 2$, we have
$$|\langle \mu_n^a-\mu,\varphi\rangle|\leq A 
\Big[1+\log^+{1\over \dist(a,\Ec)}\Big]^{\alpha/2}\|\varphi\|_{\Cc^\alpha}
\lambda^{-\alpha n/2},$$
where $A>0$ is a constant independent of $n$, $a$ and $\varphi$. 
\end{theorem}

Here, we use the notation
$$\langle\mu,\varphi\rangle:=\int\varphi d\mu.$$
The simple convergence $\mu_n^a\to\mu$ is equivalent to the convergence of the integral $\langle \mu_n^a-\mu,\varphi\rangle$  to 0. The above theorem gives us the exponential speed of convergence.

Note that the distance $\dist(a,\Ec)$ is with respect to the
Fubini-Study metric on $\P^k$. When $\Ec$ is empty (this is the case for generic $f$), by
convention, this distance is the diameter of $\P^k$ which is a finite
number. A priori, the constant $A$ depends on
$\lambda$ and $\alpha$. Note also that we have the estimate
$$A \Big[1+\log^+{1\over \dist(a,\Ec)}\Big]^{\alpha/2}
\leq A \Big[1+\log^+{1\over\dist(a,\Ec)}\Big]$$
for $0<\alpha\leq 2$. 

It is known that the measure $\mu$ has no mass on algebraic sets, in particular, on $\Ec$. So, the above result is optimal in the sense that $\mu_n^a$ does not converge to $\mu$ when $a\in\Ec$. We can show that there is a finite family of probability measures, independent of $a$, such that any limit value of $\mu_n^a$ is an element of this family. However, the choice of the limit measures depends on $a$.
For example, in dimension $k=1$, if $f(z)=z^{-d}$, we have $\Ec=\{0,\infty\}$, $\mu$ is the Haar measure on the unit circle and the above family contains three measures: $\mu$, the Dirac mass at $0$ and the Dirac mass at $\infty$.

We also deduce from the above theorem that $\mu_n^a$
converges to $\mu$ locally uniformly for $a\in\P^k\setminus\Ec$. 
The simple convergence without speed estimate was obtained in dimension 1
by Brolin \cite{Brolin} for
polynomials, by Lyubich \cite{Lyubich}, Freire-Lopes-Ma{\~n}{\'e} 
\cite{FLM} for general maps and in higher dimension by Forn\ae ss-Sibony
\cite{FornaessSibony1}, Briend-Duval 
\cite{BriendDuval1} and Dinh-Sibony \cite{DinhSibony1}.

The following corollary gives us a geometric interpretation of the above result.

\begin{corollary}
Let $U$ be an open subset of $\P^k$ such that $\mu$ has no mass on the boundary of $U$. Then, if $a$ is a point outside $\Ec$, we have
$$\# (f^{-n}(a)\cap U)=\mu(U) d^{kn}+ o(d^{kn}).$$
\end{corollary}

So, if $\mu(U)>0$, $a,b$ are two generic points and $n$ is large enough, the number of points of $f^{-n}(a)$ in $U$ is almost equal to the same quantity associated to $b$, i.e. 
$$\lim_{n\to\infty}{\#(f^{-n}(a)\cap U) \over \#(f^{-n}(b)\cap U)}=1.$$

We have the following version of Theorem \ref{th_equi_2} which is in our opinion more important.

\begin{theorem} \label{th_equi_1} 
Let $f$, $\mu$ and $\mu_n^a$ be as above. Let $1<\lambda<d$ be a fixed constant.
There is an invariant proper algebraic subset $E_\lambda$,
possibly empty, of
$\P^k$ such that if $a$ is a point out of $E_\lambda$ and if 
$\varphi$ is a $\Cc^\alpha$ function on $\P^k$
with $0<\alpha \leq 2$, then
$$|\langle \mu_n^a-\mu,\varphi\rangle|\leq A 
\Big[1 +\log^+{1\over \dist(a,E_\lambda)}\Big]^{\alpha/2}\|\varphi\|_{\Cc^\alpha}
\lambda^{-\alpha n/2},$$
where $A>0$ is a constant independent of $n$, $a$, $\varphi$.
\end{theorem}

From Theorem \ref{th_equi_1}, we can deduce several fundamental statistical properties of the measure $\mu$. 
Recall that locally this measure can be written as a Monge-Amp\`ere measure with H\"older continuous potential. It is shown by Nguyen and the authors in \cite{DinhNguyenSibony} that $\mu$ is moderate, i.e. it satisfies some exponential estimate for plurisubharmonic functions \`a la H\"ormander as in Lemma \ref{lemma_horm} where we replace $\omega_\FS^k$ by $\mu$ and $|\varphi|$ by a constant times $|\varphi|$. Therefore, in our setting, we can work with $\mu$ as if it were the Lebesgue measure. 

Theorem \ref{th_equi_1} implies a slightly weaker estimate than the following exponential mixing of $\mu$ which were proved in \cite{DinhSibony1, FornaessSibony1} : if $\varphi$ is a test $\Cc^\alpha$ function with $0<\alpha\leq 2$ and $\psi$ is a function in $L^\infty(\mu)$ then
$$\big|\langle \mu, \varphi(\psi\circ f^n)\rangle-\langle \mu,\varphi\rangle\langle\mu,\psi\rangle\big|\leq 
 A\|\varphi\|_{\Cc^\alpha}\|\psi\|_\infty d^{-\alpha n/2}.$$
 
The mixing implies the ergodicity and then, by Birkhoff's theorem, if $a$ is a $\mu$-generic point in $\P^k$, the orbit of $a$ is equidistributed in the support of $\mu$. More precisely, we have
$${1\over n}\big(\delta_a+\delta_{f(a)}+\cdots+\delta_{f^{n-1}(a)}\big)\to \mu.$$

We can also deduce more precise informations about this convergence, namely, it is possible to obtain the central limit theorem and the large deviations theorem which were proved in \cite{DinhNguyenSibony, DinhSibony4}.

\section{Sketch of the proof of Theorem \ref{th_equi_1}} \label{section_proof}

The use of Proposition \ref{prop_horm} below is new and it simplifies the original proof of Theorem \ref{th_equi_1}. Note also that Theorem \ref{th_equi_2} is a consequence of the last one. For the details, we refer to \cite{DinhSibony5}.
The main tool we use is pluripotential theory.  We recall here some results and refer to
\cite{Demailly,DinhSibony4,Hormander} for the details.

Let $\varphi:X\to \R\cup\{\infty\}$ be a function on a connected complex manifold $X$ which is not identically $-\infty$. It is called {\it plurisubharmonic} (p.s.h. for short) if its restriction to any holomorphic disc is either subharmonic or equal to $-\infty$.  It is called {\it quasi-p.s.h.} if it is locally the difference of a p.s.h. function with a smooth function. 
So, $\Cc^2$ functions are quasi-p.s.h.
A set
$E$ in $\P^k$ is {\it pluripolar} if it is contained in the pole set
$\{\varphi=-\infty\}$ of a quasi-p.s.h. function $\varphi$. 

Recall that a function $\varphi$ on $\P^k$, defined out of a
pluripolar set, is {\it d.s.h.} if it is equal to the difference
of two quasi-p.s.h. functions. We identify two d.s.h. functions if
they are equal outside of a pluripolar set. 
We summarize here some properties of these functions, see \cite{DinhSibony4} for details.
If $\varphi$ is d.s.h., 
there are two positive closed $(1,1)$-currents $S^\pm$ of the same mass such that
$\ddc\varphi= S^+-S^-$. Conversely, if $S^\pm$ are positive closed $(1,1)$-currents of the same mass, there is a d.s.h. function $\varphi$, unique up to a constant, such that $\ddc\varphi=S^+-S^-$. 

In what follows, we only consider the space $\Fc$  of d.s.h. functions $\varphi$ such that $\langle \mu,\varphi\rangle =0$, where $\mu$ is the equilibrium measure of $f$. This space is endowed with the following norm
$$\|\varphi\|:=\inf \big\{\|S^\pm\|,\quad S^\pm \mbox{ positive closed } (1,1) \mbox{-currents such that } \ddc\varphi=S^+-S^-\big\}.$$
The classical exponential estimate for p.s.h functions implies the following.

\begin{lemma} \label{lemma_horm}
There is a positive constant $C$ such that if $\varphi$ is a function in $\Fc$ with $\|\varphi\|\leq 1$, then
$$\langle \omega_\FS^k,e^{|\varphi|}\rangle \leq C.$$
\end{lemma}

The following consequence is crucial in the proof of Theorem \ref{th_equi_1}. It
is already interesting when $U=\P^k$. 
The estimate can be extended to $\varphi$ in any compact family of d.s.h. functions.
It is important to observe that we get a pointwise estimate and this allows us to get an analog in the infinite dimensional case, i.e. for super-potentials. 

\begin{proposition} \label{prop_horm}
Let $\varphi$ be a function in $\Fc$ such that $\|\varphi\|\leq 1$. Assume that $\varphi$ is H\"older
continuous on an open set $U$: $|\varphi(x)-\varphi(y)|\leq
M\dist(x,y)^\beta$ for some constants $M\geq 1$, $0<\beta\leq 1$ and
for $x,y$ in $U$.  
Then, there is a constant
$A_0>0$ independent of $\varphi$, $U$, $M$ and $\beta$ such that 
$$|\varphi(a)|\leq A_0\beta^{-1}(1+\log M)$$
for every point $a$ such that $\dist(a,\P^k\setminus U)\geq M^{-1/\beta}$.
\end{proposition}
\proof
Let $A_0>2$ be a constant large enough. If the above estimate were
false, then 
there is a function $\varphi$ and a point $a$ as above such
that 
$|\varphi(a)|\geq A_0\beta^{-1}(1+\log M)$. So, the ball $B$ of center $a$ and
of radius $M^{-1/\beta}$ is contained in $U$. We deduce from the
H\"older continuity of $\varphi$ that for every $b\in B$ 
$$|\varphi(b)|\geq A_0\beta^{-1}(1+\log M)-1\geq {1\over 2}A_0\beta^{-1}(1+\log M)\geq {1\over 2} A_0 +{1\over 2} A_0\beta^{-1}\log M.$$
This contradicts the
exponential estimate in the previous lemma.
\endproof

The mass of a positive closed $(1,1)$-current $S$ in $\P^k$ is defined by $\|S\|:=\langle S,\omega_\FS^{k-1}\rangle$. It depends only on the cohomology class of $S$ in $H^{1,1}(\P^k,\C)\simeq \C$. Using the B\'ezout theorem, we can show  that $f_*$ acts on $H^{1,1}(\P^k,\C)$ as multiplication by $d^{k-1}$. We can deduce from these properties and the total invariance of the measure $\mu$ the following lemma.

\begin{lemma} \label{lemma_lambda}
The endomorphism $f$ induces a linear operator $f_*:\Fc\to\Fc$ such that $\|f_*\|\leq d^{k-1}$. 
\end{lemma}

Recall that if $\varphi$ is a function on $\P^k$ then the function $f_*(\varphi)$ is defined by
$$f_*(\varphi)(a):=\sum_{b\in f^{-1}(a)} \varphi(b),$$
where the points in $f^{-1}(a)$ are counted with multiplicity. For $\varphi$ continuous, $f_*(\varphi)$ is continuous. 
If $\varphi$ is an $L^1$ function and $\nu$ is the Radon measure given by a smooth volume form, we have
$$\langle \nu, f_*(\varphi)\rangle = \langle f^*(\nu), \varphi\rangle.$$
For a general Radon measure $\nu$, we can define another Radon measure $f^*(\nu)$ using the same identity with
$\varphi$ continuous. 

We now define the exceptional set $E_\lambda$. 
Let $\kappa_n(x)$ denote the multiplicity of $f^n$ at $x$, i.e. the local topological degree of $f^n$ at $x$,
for $n\geq 0$. More precisely, for $z$ generic near $f^n(x)$,
$f^{-n}(z)$ has $\kappa_n(x)$ points near $x$. 
Define 
$$\kappa_{-n}(x):=\max_{y\in f^{-n}(x)}\kappa_n(y).$$
It was shown in \cite{Dinh} that the sequence
$\kappa_{-n}^{1/n}$ converges to a
function $\kappa_-$ which is upper semi-continuous with respect to the
Zariski topology. 

Moreover, for any $\delta>1$, the level set
$\{\kappa_-\geq \delta\}$ is an invariant proper algebraic subset of
$\P^k$. Define $E_\lambda:=\{\kappa_-\geq d/\lambda\}$. So, there is a constant $1<\delta_0<d/\lambda$ such that $\kappa_{n_0} < \delta_0^{n_0}$ outside $f^{-n_0}(E_\lambda)$ for a fixed integer $N_0$ large enough. 
In what follows, without loss of generality, we replace $f,d,\lambda,\delta_0$ by $f^n,d^n,\lambda^n,\delta_0^n$ for $n$ large enough in order to assume  that the multiplicity of $f$ at any point outside $f^{-1}(E_\lambda)$ is smaller than $\delta_0$ and that $20k^2\delta_0<d/\lambda$. 

The previous property of multiplicity allows us to prove a version of the classical Lojasiewicz inequality adapted to our situation. Denote by  $V_t$ the $t$-neighbourhood of $E_\lambda$.

\begin{proposition} \label{prop_loj}
There is an integer $N\geq 1$ and a constant $A_1\geq 1$ such that if $0<t<1$ is a
constant and if $x,y$ are two points outside $V_t$, then we can write 
$$f^{-1}(x)=\{x_1,\ldots, x_{d^k}\}\quad \mbox{and}\quad 
f^{-1}(y)=\{y_1,\ldots, y_{d^k}\}$$
with $\dist(x_i,y_i)\leq A_1t^{-N}\dist(x,y)^{1/\delta_0}$.
\end{proposition}

Finally, we have the following proposition where we assume that $\varphi$ is a function in $\Fc$ such that 
$\|\varphi\|_{\Cc^2}\leq 1$. Define $\Lambda:=d^{1-k}f_*$.

\begin{proposition}\label{prop_holder_Lambda}
The function $\Lambda^n(\varphi)$ is H\"older continuous on
$\P^k\setminus E_\lambda$. Moreover, there is a constant $A_2\geq 1$
such that for every $n\geq 0$ and $0<t\leq 1$
$$|\Lambda^n\varphi(x)-\Lambda^n\varphi(y)|\leq
A_2^{Nn^2}t^{-Nn}\dist(x,y)^{\delta_0^{-n}}$$
for $x,y$ out of $V_t$.
\end{proposition}
\proof
Since $E_\lambda$ is invariant, it is not difficult to see that for a constant $c>0$ small enough, we have
$$\dist(f^{-1}(x), E_\lambda)\geq c \dist(x,E_\lambda).$$
Define $A_2:=dA_1/c$.  
The proof is by induction on $n$. The case $n=0$ is clear.
Assume that the proposition holds for $n$. We show it for $n+1$.
Let $x,y$ be two points out of $V_t$. By Proposition \ref{prop_loj}, we can write 
$$f^{-1}(x)=\{x_1,\ldots,x_{d^k}\}\quad \mbox{and}\quad 
f^{-1}(y)=\{y_1,\ldots,y_{d^k}\}$$
so that $\dist(x_i,y_i)\leq A_1t^{-N}\dist(x,y)^{1/\delta_0}$. Observe that $x_i$ and $y_i$ are out of 
$V_{ct}$. We deduce from the definition of $\Lambda$ and the
induction hypothesis that 
\begin{eqnarray*}
|\Lambda^{n+1}\varphi(x)-\Lambda^{n+1}\varphi(y)| & \leq & d^{1-k}\sum
|\Lambda^n\varphi(x_i) -\Lambda^n\varphi(y_i)| \\
&\leq& d^{1-k}A_2^{Nn^2}(ct)^{-Nn}\sum
\dist(x_i,y_i)^{\delta_0^{-n}} \\
& \leq & dA_2^{Nn^2} c^{-Nn}
t^{-Nn}A_1^{\delta_0^{-n}}t^{-N\delta_0^{-n}}\dist(x,y)^{\delta_0^{-n-1}}\\
&\leq & A_2^{N(n+1)^2}t^{-N(n+1)}\dist(x,y)^{\delta_0^{-n-1}}.
\end{eqnarray*}
This completes the proof.
\endproof

\noindent
{\bf End of the proof of Theorem \ref{th_equi_1}.}
First observe that from the theory of interpolation between Banach spaces (in our case between $\Cc^0$ and $\Cc^2$), it is enough to consider the case $\alpha=2$. Indeed, if $L$ is a continuous linear form on the space of continuous functions on $\P^k$, it defines also a continuous linear form on $\Cc^\alpha$ and we have the inequality
$$\|L\|_{\Cc^\alpha}\leq A\|L\|_\infty^{1-\alpha/2} \|L\|_{\Cc^2}^{\alpha/2},$$
where $A>0$ is a constant independent of $L$, see \cite{Triebel} for details. In our situation, it is enough to apply this inequality to the Radon measures $d^{-kn}(f^n)^*(\delta_a)-\mu$.

So, assume that $\alpha=2$ and $\varphi$ is a function of class $\Cc^2$. 
Since the theorem is clear for constant test functions, subtracting from $\varphi$ a constant allows us to assume that $\varphi$ is a function in $\Fc$. Moreover, by linearity, we can assume that $\|\varphi\|\leq 1$ and $\|\varphi\|_{\Cc^2}$ is bounded. We use here that $\|\ \|\lesssim \|\ \|_{\Cc^2}$. 
By Lemma \ref{lemma_lambda}, we have $\|\Lambda^n(\varphi)\|\leq 1$.

We also obtain from the definition of $\Lambda$ that 
$$\langle d^{-kn} (f^n)^*(\delta_a),\varphi\rangle=d^{-n}\Lambda^n\varphi(a).$$
Define 
$$l:=1+\log^+{1\over \dist(a,E_\lambda)}\cdot$$ 
We need to show that
$|\Lambda^n\varphi (a)|\leq A l\lambda^{-n}d^n$
for some constant $A>0$ and for $n\geq 1$. Define $t:=e^{-l}$. Observe
that 
$\dist(a,V_t) \geq t$. 
Therefore, Propositions \ref{prop_holder_Lambda} and  \ref{prop_horm} yields
$$|\Lambda^n\varphi(a)|\leq
A_0\delta_0^n\big[1+\log(A_2^{Nn^2}t^{-Nn})\big]\lesssim
l\lambda^{-n}d^n$$
since $\delta_0<d/\lambda$ and $l\geq 1$.  
\hfill $\square$

\section{Equidistribution of varieties} \label{section_varieties}

In this section we survey the results on equidistribution of varieties. Recall that the sequence $d^{-n}(f^n)^*(\omega_\FS)$ converges to a canonical invariant positive closed $(1,1)$-current $T$ of mass 1. We call it {\it the Green current} of $f$. For any integer $1\leq p\leq k$, the sequence $d^{-pn}(f^n)^*(\omega_\FS^p)$ converges to the positive closed $(p,p)$-current $T^p:=T\wedge\ldots\wedge T$ ($p$ times) that we call {\it the Green current of order $p$} or {\it the Green $(p,p)$-current} of $f$, see e.g. \cite{DinhSibony4}. When $p=k$ we obtain the equilibrium measure $\mu=T^k$ considered above.  

Note that the operators $f_*$ and $f^*$ are well-defined on positive closed currents \cite{DinhSibony6}. If $S$ is a positive closed $(p,p)$-current on $\P^k$, we have $\|f^*(S)\|=d^p\|S\|$ and $\|f_*(S)\|=d^{k-p}\|S\|$.  If $V$ is an algebraic set of pure codimension $p$, the integration on its regular part defines a positive closed $(p,p)$-current $[V]$. The mass of $[V]$ is equal to the degree of $V$. We conjecture the following.

\begin{conjecture} \label{conjecture_varieties}
Let $V$ be generic in the Zariski sense among algebraic sets of pure codimension $p$ and of a given degree. 
Then, $\deg(V)^{-1}d^{-pn}(f^n)^*[V]$ converge to $T^p$ in the sense of currents when $n$ goes to infinity. Moreover, the convergence is exponentially fast. 
\end{conjecture}

Fix a constant $1<\lambda<d$. We expect a much stronger property: there is a finite family of algebraic sets $E_\lambda^1,\ldots,E_\lambda^m$ such that if $V$ intersects $E_\lambda^i$ properly, i.e. the intersection is empty when $\dim E_\lambda^i<p$ and the intersection is of dimension $\dim E_\lambda^i-p$ otherwise, then
$$\big|\langle \deg(V)^{-1}d^{-pn} (f^n)^*[V]-T^p,\Phi\rangle \big|\lesssim \|\Phi\|_{\Cc^\alpha} \lambda^{-n\alpha/2}$$
for any test form $\Phi$ of class $\Cc^\alpha$ with $0<\alpha\leq 2$.

If we replace the family of $E_\lambda^i$ by the family of all totally invariant algebraic sets, we also expect the exponential convergence with rate $\lambda^{-n\alpha/2}$ for some $\lambda>1$. 

The conjecture says in particular that if $U$ is an open set such that $T^p$ has positive mass on $U$ but no mass on $\partial U$ and if $V,V'$ are two generic algebraic sets of codimension $p$, of the same degree, then the volume of $f^{-n}(V)\cap U$ is almost equal to the one of $f^{-n}(V')\cap U$ when $n$ is large enough. Here, by volume we mean the Hausdorff $2(k-p)$-dimensional measure with respect to a fixed Hermitian metric on $\P^k$.

We have seen that the conjecture holds for the case of points, i.e. $p=k$. The conjecture was recently confirmed in the case of hypersurfaces, i.e. $p=1$, by Taflin in \cite{Taflin}. Taflin's theorem is as follows.

\begin{theorem}
Let $\lambda$ be a constant such that $1<\lambda<d$. There is a finite family of algebraic sets $E_\lambda^1,\ldots,E_\lambda^m$ satisfying the following property. If $V$ is a hypersurface 
which does not contain any $E_\lambda^i$, then
$$\big|\langle \deg(V)^{-1}d^{-n} (f^n)^*[V]-T,\Phi\rangle \big|\leq C \|\Phi\|_{\Cc^\alpha} \lambda^{-n\alpha/2},$$
for every test form $\Phi$ of class $\Cc^\alpha$ with $0<\alpha\leq 2$, where $C>0$ is a constant depending on $f,\lambda,V$ and $\alpha$. 
\end{theorem}

We refer to the original paper by Taflin for  the explicit construction of $E_\lambda^i$. Note that the convergence without speed was proved by Forn\ae ss-Sibony \cite{FornaessSibony2} and Favre-Jonsson \cite{FavreJonsson} in dimension $k=2$ and by Dinh-Sibony \cite{DinhSibony2} in any dimension. In this case, we can replace the  family of $E^i_\lambda$ by the family of minimal totally invariant algebraic sets, see also Para \cite{Para}. 

A key point in the proof is to write, in a unique way, 
$$\deg(V)^{-1}[V]-T=\ddc u$$
with $u$ a function in $\Fc$. The uniqueness of the solution of such an equation and the total invariance of $T$ imply that
$$\deg(V)^{-1}d^{-n} (f^n)^*[V]-T=\ddc (d^{-n}u\circ f^n).$$
So, the problem is reduced to prove the convergence of the sequence of functions $d^{-n}u\circ f^n$ to 0 in a weak sense and to estimate the speed of convergence.

Taflin proved and used some version of exponential estimates for non-compact families of d.s.h. functions. Moreover, the proof contains an induction on the dimension, that is, he has to check the convergence on some algebraic subsets of $\P^k$ which may be singular. This requires precise versions of Lojasiewicz inequality in order to handle different technical difficulties.

The above conjecture is still open in the general case. We have nevertheless the following result \cite{DinhSibony3} which proves the conjecture for generic maps\footnote{The speed of convergence is not precisely stated in \cite{DinhSibony3}, but the proof gives Theorem \ref{th_equi_v}.}. Recall that the parameter space for holomorphic endomorphisms of algebraic degree $d$ is a connected complex 
quasi-projective manifold $\Hc_d(\P^k)$.

\begin{theorem} \label{th_equi_v}
Let $1<\lambda<d$ be a constant. There is a non-empty Zariski open set $\Hc_d^\lambda(\P^k)$ of $\Hc_d(\P^k)$ such that if $f$ is an element in this open set, then  $\deg(V)^{-1}d^{-pn}(f^n)^*[V]$ converge to $T^p$ in the sense of currents for any algebraic set $V$ of pure codimension $p$. Moreover, we have 
$$\big|\langle \deg(V)^{-1}d^{-pn} (f^n)^*[V]-T^p,\Phi\rangle \big|\leq C \|\Phi\|_{\Cc^\alpha} \lambda^{-n\alpha/2},$$
for every test form $\Phi$ of class $\Cc^\alpha$ with $0<\alpha\leq 2$. Here, $C$ is a positive constant independent of $V$.   
\end{theorem}

The set $\Hc_d^\lambda(\P^k)$ is defined as a set of maps where the local multiplicity is not too big.

A main difficulty is that the above theory of plurisubharmonic functions is not enough to handle algebraic cycles of arbitrary dimension. So, for this purpose, we introduced and developed a theory of super-potentials for positive closed currents and we applied it in the dynamical setting. Roughly speaking, this is a theory of quasi-p.s.h. functions in infinite dimension. We are able for example to get a version of the exponential estimate in Proposition \ref{prop_horm}. In order to obtain the solution of Conjecture \ref{conjecture_varieties}, we still have to get a version of the Lojasiewicz inequality which seems to be a difficult problem.


T.-C. Dinh, UPMC Univ Paris 06, UMR 7586, Institut de
Math{\'e}matiques de Jussieu, F-75005 Paris, France. {\tt
  dinh@math.jussieu.fr}, {\tt http://www.math.jussieu.fr/$\sim$dinh} 

\

\noindent
N. Sibony,
Universit{\'e} Paris-Sud, Math{\'e}matique - B{\^a}timent 425, 91405
Orsay, France. {\tt nessim.sibony@math.u-psud.fr}


\small
\begin{thebibliography}{99}   


\bibitem{AmerikCampana}
Amerik E. and Campana F.,
Exceptional points of an endomorphism of the projective plane,
{\it Math. Z.}, {\bf 249} (2005), no. 4, 741-754. 
 
\bibitem{BriendDuval1}
Briend J.-Y. and Duval J., Deux caract{\'e}risations de la mesure
d'{\'e}quilibre d'un endomorphisme de ${\rm P}\sp k(\bold C)$, {\it
Publ. Math. Inst. Hautes {\'E}tudes Sci.}, {\bf 93} (2001), 145-159.


\bibitem{Brolin}
Brolin H., Invariant sets under iteration of rational functions,
{\it Ark. Mat.}, {\bf 6} (1965), 103-144.


\bibitem{CerveauLinsNeto}
Cerveau D. and Lins Neto A.,
Hypersurfaces exceptionnelles des endomorphismes de ${\rm CP}(n)$,
{\it Bol. Soc. Brasil. Mat. (N.S.)}, {\bf 31} (2000), no. 2, 155-161. 

\bibitem{Demailly}
Demailly J.-P., 
{\it Complex analytic and differential geometry}, available at \\
{\tt www.fourier.ujf-grenoble.fr/$\sim$demailly}. 


\bibitem{Dinh}
Dinh T.-C., Analytic multiplicative cocycles over holomorphic dynamical
systems, {\it special issue of 
Complex Variables and Elliptic Equations}, {\bf 54} (2009), no. 3-4, 243-251.


\bibitem{DinhNguyenSibony}
Dinh T.-C., Nguyen V.-A. and Sibony N., 
Exponential estimates for plurisubharmonic functions and
stochastic dynamics, {\it J. Differential Geom.}, {\bf 84} (2010), no. 3, 465-488. 
 
 

\bibitem{DinhSibony1}
Dinh T.-C. and Sibony N., Dynamique des applications d'allure
polynomiale, {\it J. Math. Pures  Appl.}, {\bf 82} (2003),
367-423.

\bibitem{DinhSibony6}
------, Pull-back of currents by holomorphic maps, {\it
  Manuscripta Math.,} {\bf 123} (2007), 357-371.

\bibitem{DinhSibony2}
------, Equidistribution towards the Green current for
holomorphic maps, {\it Ann. Sci. {\'E}cole Norm. Sup.}, {\bf 41} (2008),
307-336.


\bibitem{DinhSibony3}
------, Super-potentials of positive closed currents, intersection theory
and dynamics, {\it Acta Math.}, {\bf 203} (2009), no. 1, 1-82.


\bibitem{DinhSibony5}
------, Equidistribution speed for endomorphisms of projective spaces, {\it Math. Ann.}, {\bf 347} (2010), no. 3, 613-626.

\bibitem{DinhSibony4}
------, Dynamics in several complex variables: endomorphisms of
  projective spaces and polynomial-like mappings, 165-294, {\it Lecture Notes in Math.}, {\bf 1998}, Springer, Berlin, 2010.

\bibitem{FavreJonsson}
Favre C. and Jonsson M., Brolin's theorem for curves in two complex
dimensions,  {\it Ann. Inst. Fourier (Grenoble)}, {\bf  53}  (2003),  no. 5, 1461-1501.

\bibitem{Fornaess}
Forn\ae ss J.-E., Dynamics in several complex variables, {\it CBMS Regional Conference Series in Mathematics}, 
{\bf 87},  American Mathematical Society, Providence, RI, 1996. 

\bibitem{FornaessSibony1}
Forn\ae ss J.-E. and Sibony N., 
Complex dynamics in higher dimensions. Notes partially written by
Estela A. Gavosto, 
\textit{NATO Adv. Sci. Inst. Ser. C Math. Phys. Sci.}, {\bf 439},  Complex potential
theory (Montreal, PQ, 1993),  131-186, {\it Kluwer Acad. Publ.}, Dordrecht, 1994.

\bibitem{FornaessSibony3}
------, Complex dynamics in higher dimension. I. Complex analytic methods in dynamical systems (Rio de Janeiro, 1992), {\it Ast\'erisque}, {\bf 222} (1994), 5, 201-231.

\bibitem{FornaessSibony2}
------, 
Complex dynamics in higher dimension.
II.  Modern methods in complex analysis (Princeton, NJ, 1992),
135-182, {\it Ann. of Math. Stud.}, {\bf 137}, Princeton Univ.
Press, Princeton, NJ, 1995.

\bibitem{FLM}
Freire A., Lopes A. and Ma{\~n}{\'e} R., An invariant measure for
rational maps, {\it Bol. Soc. Brasil. Mat.}, {\bf 14} (1983), no.
1, 45-62.


\bibitem{Hormander}
H{\"o}rmander L., {\it An introduction to complex analysis in several
  variables}, Third edition, 
North-Holland Mathematical Library, {\bf 7}, North-Holland Publishing
Co., Amsterdam, 1990.


\bibitem{Lyubich}
Lyubich M. Ju., Entropy properties of rational endomorphisms of
the Riemann sphere, {\it Ergodic Theory Dynam. Systems}, {\bf  3}
(1983), no. 3, 351-385.

\bibitem{Para}
Para M.R., The Jacobian cocycle and equidistribution towards the Green current, {\it preprint}, 2011. 
{\tt 	arXiv:1103.4633v1}

\bibitem{Sibony}
Sibony N., Dynamique des applications rationnelles de
$\mathbb{P}^k$, \textit{Panoramas et Synth{\`e}ses},  \textbf{8}
(1999), 97-185.

\bibitem{Taflin}
Taflin J., Equidistribution speed towards the Green current for endomorphisms of $\P^k$, {\it Advances in Math.}, to appear. {\tt arXiv:1011.0641}

\bibitem{Triebel}
Triebel H.,
{\it Interpolation theory, function spaces, differential operators},
North-Holland, 1978. 

\bibitem{Yuan}
Yuan X., Algebraic Dynamics, Canonical Heights and Arakelov 
Geometry, {\it preprint}, 2011.

\end{thebibliography}
\end{document}